\documentclass[a4,11pt]{amsart}
\usepackage{latexsym}
\usepackage[dvips]{graphicx,color,epsfig}
\usepackage{lscape, rotating}
\usepackage{varioref,float,amssymb,amsmath}
\usepackage{float,amsfonts,amsthm}
\usepackage{changebar,longtable}
\usepackage{pstricks,pstricks-add,pst-math,pst-xkey}
\usepackage{psfrag} 
\usepackage{overpic} 
\newtheorem{obs}{Remark}[section]

\newtheorem{prop}[obs]{Proposition}
\newtheorem{teo}[obs]{Theorem}
\newtheorem{lem}[obs]{Lemma}

\def\dem{{\bf Proof:}\hspace{.2in}}
\begin{document}
\author{
J. Bastos, T. Rodrigues de Souza}

\address{Departamento de MatemÃÂ¡tica, UNESP - Universidade Estadual Paulista, Rua Crist\'ov\~ao Colombo, 2265, Jardim Nazareth, 15054-000, S\~ao Jos\'e do Rio Preto, SP, Brazil.}
\email{jeferson@ibilce.unesp.br}


\address{Departamento de Matem\'atica, UNESP - Universidade Estadual Paulista, AV. Eng. Luiz Ed. Carrijo Coube, 14-01, Vargem Limpa, 17033-360, Bauru, SP, Brazil.}
\email{tatimi@fc.unesp.br}

\date{\today}
\title{Parametrization for a class of Rauzy Fractals}
\begin{abstract}
In this paper, we study a class of Rauzy fractals ${\mathcal R}_a$
given by the polynomial $x^3- ax^2+x-1$ where $a \geq 2$ is an
integer. In particular, we give explicitly an automaton that
generates the boundary of ${\mathcal R}_a$ and using an unusual
numeration system we prove that ${\mathcal R}_a$ is homeomorphic to
a topological disk.
\end{abstract}
\subjclass[2000]{} 


\maketitle

\section{Introduction}
The Rauzy fractal is a compact subset of $\mathbb{R}^{n},\ n\geq 1$.
It was studied by many mathematicians and is connected to many
topics such as: numeration systems (\cite{BBLT},\cite{BeSi},
\cite{b},\cite{m1}), geometrical representation of symbolic
dynamical systems (\cite{a1}, \cite{m2}), multidimensional continued
fractions and simultaneous approximations (\cite{a2}, \cite{c}),
auto-similar tilings (\cite{a1}, \cite{b}), substitutions and
tilings (\cite{BST}) and Markov partitions of Hyperbolic
automorphisms of Torus (\cite{m2}, \cite{b}).




    Let $\beta>1$ be a fixed real number. Any positive real number $x$ can be expanded as $$x=\sum_{i=N_0}^{\infty}a_{-i}\beta{-i}=a_{-N_0}\beta^{-N_0}+a_{-N_0-1}\beta^{-N_0-1}+...$$ with $a_{i}\in\mathbb{Z}\cap[0,\beta)$ and we are assuming the greedy condition $$ \left|x-\sum_{i=N_0}^{N}a_{-i}\beta^{-i}\right|<\beta^{-N},$$ for all $N\geq N_0$. We call this expansion a beta expansion of $x$ in base $\beta$. A Pisot number is an algebraic integer whose conjugates other than itiself have modulus less than one. Let $Fin(\beta)$ be a set consisting of all finite beta expansions and consider the condition $$Fin(\beta)=\mathbb{Z}[\beta^{-1}]_{\geq 0}. \ \ \   (F)$$
    Consider the beta expansion of the positive number $$0<1-[\beta]\beta^{-1}=c_{-2}\beta^{-2}+c_{-3}\beta^{-3}+...=.0c_{-2}c_{-3}...$$
    If we put $c_{-1}=[\beta]$ we can write $$1=.c_{-1}c_{-2}c_{-3}...$$ This expansion $.c_{-1}c_{-2}c_{-3}...$ is called the expansion of $1$ and denoted by $d(1,\beta)$. We can identify this expression with the word $c_{-1}c_{-2}c_{-3}...$ generated by $\mathbb{A}=\mathbb{Z}\cap[0,\beta)$. Every finite word generated by $\mathbb{A}$ represents a beta expansion in base $\beta$ if and only if the word is lexicographically less than $d(1,\beta)$ at any starting point. This fact can be generalized to infinite words apart from certain exceptions (see \cite{Parry}).

        In \cite{FS92} they proved that if $\beta>1$ is an integer then $(F)$ holds and, conversely, the condition $(F)$ implies that $\beta$ is a Pisot number. A Pisot number $\beta$ is called a Pisot unit if it is also a unit of the integer ring of $\mathbb{Q}[\beta]$.  In \cite{akiyama1} we have the following results :
        \begin{teo}\label{th1.1}
            Let $\beta$ ba a cubic Pisot number. Then $\beta>1$ has property $(F)$ if and only if $\beta$ is a root of the following polynomial with integer coefficients:$$x^3-ax^2-bx-1,a\geq 0,\mbox{ and }-1\leq b \leq a+1.$$
        \end{teo}

            \begin{lem}\label{lem1.2}
                Let $\beta>1$ be a cubic Pisot number with $Irr(\beta)=x^3-ax^2-bx-1$. Then the expansion of $1$ in base $\beta$ is given by:
                \begin{itemize}
                    \item[$i)$] $d(1,\beta)=.(a-1)(a+b-1)\widetilde{(a+b)}$, if $-a+1\leq b \leq -2;$
                    \item[$ii)$] $d(1,\beta)=.ab1$, if $0\leq b \leq a;$
                    \item[$iii)$] $d(1,\beta)=.(a-1)(a-1)01$, if $b=-1;$
                    \item[$iv)$] $d(1,\beta)=.(a+1)00a1$, if $b=a+1.$
                \end{itemize}
                Here $\widetilde{w}$ is the periodic expansion $wwww...$.
            \end{lem}
            \begin{teo}
                A cubic Pisot unit $\beta$ has property $(F)$ if and only if $d(1,\beta)$ is finite.
            \end{teo}


To each cubic Pisot unit  satisfiyng $(F)$, we can associate a Rauzy
fractal. In the case where $\beta$ is a Pisot number satisfying
condition $(ii)$ Lemma \ref{lem1.2}, the Rauzy fractal was studied
in \cite{TLMS} and \cite{Loridant}. In \cite{TLMS} the authors
proved that if $2b>a-3,$ then the boundary of the Rauzy fractal is
not homeomorfic to a circle. If $\beta$ is a cubic Pisot unit
satisfiyng $ (iii)$ of Lemma  \ref{lem1.2} and $\alpha,
\overline{\alpha}$ its Galois conjugates the fractal associated is
given by
$$
\mathcal{R}_a=\left\{\sum_{i=2}^{\infty} a_{i} \alpha^{i},\;
a_{i}a_{i-1}a_{i-2}a_{i-3}<_{lex} (a-1)(a-1)01 ,\; \forall i \geq 5
\right\},$$ where $<_{lex}$ is the lexicographic order on finite
words. In \cite{b1} the authors proved the topological and
arithmetical properties of $\mathcal{R}_a$. In particular, they
proved there exists an explicit finite state automaton $\mathcal{A}$
such that the boundary of $\mathcal{R}_a$ is recognized by
$\mathcal{A}$. With this automaton they proved that for $a=2$, the
boundary of $\mathcal{R}_2$ is homeomorfic to a circle. Their proof
cannot be extended to the case $a\geq 3$. The parametrization of the
boundary of $\mathcal{R}_a,\  a\geq3$ is different from the case
$a=2$. It uses an unusual numeration system.

In this paper we will study the fractal associated to a number
$\beta$ satisfying the condition $(iii)$ of Lemma \ref{lem1.2} with
$a\geq 3, b=-1$. In this case the polinomial
 $p(x)=x^3-ax^2+x-1=(x-\beta)(x-\alpha)(x-\overline{\alpha})$,where $\beta >1$ and $\alpha,\overline{\alpha}\in\mathbb{C}\backslash\mathbb{R}$.

The purpose of this work is to present a complete description of the
boundary of $\mathcal{R}_a,\ a\geq3$. Our main result is the
following:.
\begin{teo}$\label{princtheo}\partial \mathcal{R}_{a}$ is homeomorfic to
    $S^1$.
\end{teo}
\newpage
\begin{figure}[h!]

\centering \epsfig{file=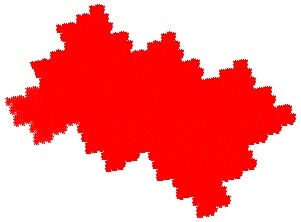,height=5cm,width=5.5cm} \caption{
{\small $\mathcal{R}_3$ }} \label{Figure2}
\end{figure}

\section{Background, notations and definitions}

In this section we will give more informations about
$\beta$-numeration, Rauzy fractal, automaton and we will present
some notations that will be used in the next sections.

Assume that $\beta$ is a Pisot number of degree $d\geq 3$. We denote
by $\beta_2,\beta_3,...,\beta_r$  the real Galois conjugate of
$\beta$ and by $\beta_{r+1},...,
\beta_{r+s},\beta_{r+s+1}=\overline{\beta_{r+1}},...,\beta_{r+2s}=\overline{\beta_{r+s}}$
its complex Galois conjugates. Let
$$\psi=(\beta_2,...,\beta_r,\beta_{r+1},...,\beta_{r+s})\in\mathbb{R}^{r-1}\times
\mathbb{C}^{s}$$ and put $\psi^i
=(\beta^i_2,...,\beta^i_r,\beta^i_{r+1},...,\beta^i_{r+s}),\ \forall
i\in\mathbb{Z}.$ The Rauzy fractal is by definition the set
$$\mathcal{R}=\left\{\sum_{i=0}^{\infty}a_i\psi^i,\ (a_i)_{i\geq
0}\in E_{\beta}\right\}$$ where $E_{\beta}=\{(x_i)_{i\geq
k},k\in\mathbb{Z}|\forall n\geq k,(x_i)_{k\leq i\leq n}\mbox{ is a
finite }\beta\mbox{ expansion}  \}$.

An important class of Pisot numbers are those such that the
associated Rauzy fractal has $0$ as an interior point. This numbers
where characterized by Akyiama (\cite{akiyama2}), and they are
exactly the Pisot numbers satisfying condition $(F)$.

In this paper we will work with sequences $(a_n)_{n\in\mathbb{Z}}$
belonging to $\{0,1,...,a-1\}^{\mathbb{Z}}$ and the following set
\\ $\mathcal{N}= \{(a_{n})_{n \in \mathbb{Z}}, \; \exists k \in
\mathbb{Z},\; a_k
>0,\;  a_i = 0 \mbox { for all } i <k,\;  a_{i}a_{i-1}a_{i-2}
a_{i-3} <_{lex} (a-1)(a-1)01,\; \forall i \geq k  \}. $

If $(a_n)\in\mathcal{N}$ we will call it an admissible sequence.

Take $(a_n,b_n)_{n\in\mathbb{Z}}$ an infinite path on the automaton
$\mathcal{A}$ starting in the initial state. If
$(a_n),(b_n)\in\mathcal{N}$  we will call it an admissible path.


%
In \cite{b1} the following results were proved.
\begin{enumerate}
\item Let $\displaystyle z=\sum_{i=2}^{\infty}a_i\alpha^{i}\in
\mathcal{R}_a$. Then $z\in\partial\mathcal{R}_a$ if and only if
there exists $(b_i)_{i\geq l}\in \mathcal{N},\ l<1,\ b_l\neq0$ such
that
$\displaystyle\sum_{i=2}^{\infty}a_i\alpha^{i}=\sum_{i=l}^{\infty}b_i\alpha^{i}.$
\item There exists an explicit finite state automaton $\mathcal{A}$
(see figure 2 below) such that
$\displaystyle\sum_{i=l}^{\infty}\epsilon_i\alpha^{i}=\sum_{i=l}^{\infty}\epsilon'_i\alpha^{i},$
$(\epsilon_i),(\epsilon'_i)\in \mathcal{N}$ if and only if
$(\epsilon_i,\epsilon'_i)_{i\geq l}$ is an admissible path.
\end{enumerate}


Let us explain the behavior of this automaton. Let $\varepsilon =
(\varepsilon_{i})_{i\geq l}$ and $\varepsilon' =
(\varepsilon'_{i})_{i\geq l}$ belonging to $\mathcal{N}$,
$x=\sum_{i=l}^{\infty}\varepsilon_{i}\alpha^{i}$ and
$y=\sum_{i=l}^{\infty}\varepsilon'_{i}\alpha^{i}$. Suppose $x=y$.
For all $k\geq l$ we put

\begin{align}
\label{def-Ak} A_k (\varepsilon ,\varepsilon') =
\alpha^{-k+2}\sum_{i=l}^{k}(\varepsilon_{i}-\varepsilon'_{i})\alpha^{i}.
\end{align}

In \cite{b1} the authors proved that $A_k (\varepsilon
,\varepsilon')\in S=\{0,\pm \alpha,\pm \alpha^2, \pm (\alpha-
\alpha^2),\;\pm (1+(a-1)\alpha^2), \pm (1+ (a-2)\alpha^2),\; \pm (1
-\alpha+ (a-1)\alpha^2),\;\pm (1-2 \alpha +a\alpha^2)\}.$ We can see
that for all $k\geq l$,

\begin{equation}\label{relaq}
A_{k+1} (\varepsilon ,\varepsilon') = \frac {A_{k}(\varepsilon
,\varepsilon')}{\alpha} + (\varepsilon_{k+1}- \varepsilon'_{k+1})
\alpha^{2} .
\end{equation}

Let $s$ be the smallest integer such that $\varepsilon_{s} \ne
\varepsilon'_{s}$. Hence $A_{i} (\varepsilon ,\varepsilon') = 0$ for
$i \in \{l, \ldots, s-1\}$. Suppose  $\varepsilon_{s}
>\varepsilon'_{s}$. Then, $A_{s}= (\varepsilon_{s}
-\varepsilon'_{s})\alpha^{2}= \alpha^{2}$. From (\ref{relaq}) we
deduce $ A_{s+1}(\varepsilon ,\varepsilon') = \alpha +
(\varepsilon_{s+1}- \varepsilon'_{s+1}) \alpha^{2}$ which should
belong to $S$. Hence $A_{s+1} (\varepsilon ,\varepsilon')= \alpha $
if $\varepsilon _{s+1}= \varepsilon'_{s+1}$ or  $A_{s+1}
(\varepsilon ,\varepsilon')= \alpha- \alpha^{2} $ if $(\varepsilon
_{s+1}, \varepsilon'_{s+1})= (t_1, t_1 +1)$, where $ 0 \leq t_1 \leq
a-2$. Continuing by the same way and using the fact that
 the set of states   $S$ is finite, we obtain the following
finite state automaton shown in Figure 2.

\begin{figure}[h!]
\begin{center}
\begin{overpic}[height=12cm]{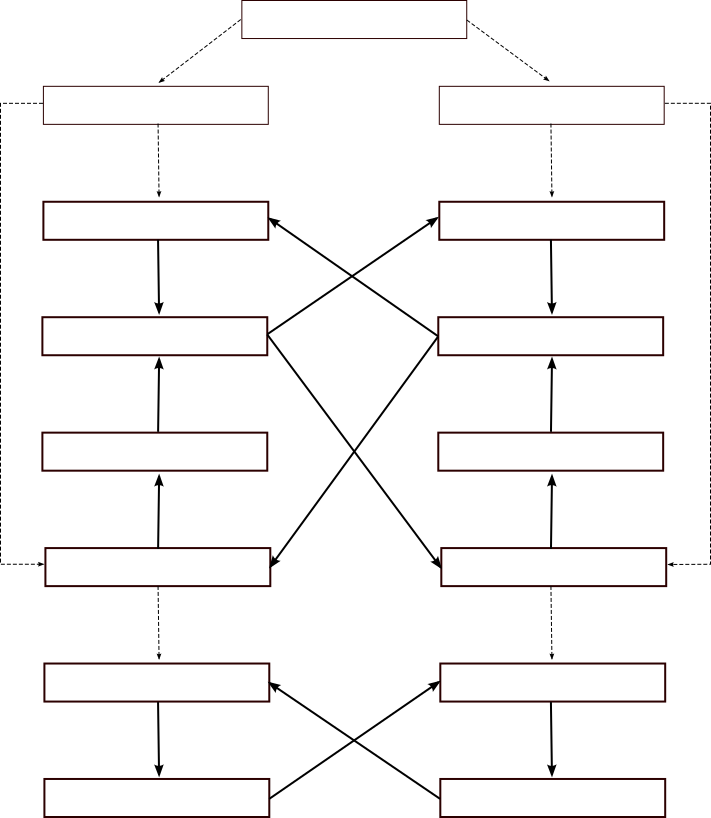} 
{\scriptsize \put(43,97){$0$}
\put(12,93){$(\varepsilon,\varepsilon+1)$}
\put(69,80){$(\varepsilon,\varepsilon+1)$} \put(68,65){$(a-1,0)$}
\put(68,51){$(\varepsilon+a-2,\varepsilon)$} \put(68,37){$(a-1,0)$}
\put(68,23){$(\varepsilon+a-2,\varepsilon)$} \put(68,8){$(a-1,0)$}
\put(8,8){$(0,a-1)$} \put(3,23){$(\varepsilon,\varepsilon+a-2)$}
\put(8,38){$(0,a-1)$} \put(3,51){$(\varepsilon,\varepsilon+a-2)$}
\put(47,74){$(\varepsilon,\varepsilon)$}
\put(33,74){$(\varepsilon,\varepsilon)$}
\put(7,79){$(\varepsilon+1,\varepsilon)$} \put(7,65){$(0,a-1)$}
\put(66,93){$(\varepsilon+1,\varepsilon)$}
\put(87,58){$(\varepsilon,\varepsilon)$}
\put(-11,58){$(\varepsilon,\varepsilon)$} \put(15,86){$-\alpha^{2}$}
\put(13,72){$-\alpha+\alpha^{2}$} \put(12,44){$-1-(a-1)\alpha^{2}$}
\put(16,30){$-\alpha$} \put(11,16){$-1-(a-2)\alpha^{2}$}
\put(11,2){$-1+2\alpha-a\alpha^{2}$}

\put(66.6,86){$\alpha^2$} \put(64,72){$\alpha-\alpha^2$}
\put(57,58){$1-\alpha+(a-1)\alpha^2$}
\put(60,44){$1+(a-1)\alpha^{2}$}
\put(51,37){$(\varepsilon+1,\varepsilon)$}
\put(26,37){$(\varepsilon,\varepsilon+1)$} \put(66,30){$\alpha$}
\put(61,16){$1+(a-2)\alpha^{2}$} \put(42,17){$(a-1,0)$}
\put(41,1){$(0,a-1)$} \put(61,2){$1-2\alpha+a\alpha^{2}$}
\put(7,58){$-1+\alpha-(a-1)\alpha^2$}}

\end{overpic}
\vspace{.5cm} \caption{Automaton $\mathcal{A}$}
\end{center}
\end{figure}

Another result proved in \cite{b1} is the following.
\begin{prop}\label{prop35}
${\mathcal R}_{a}$ induces a periodic tiling of the plane
$\mathbb{C}$ modulo $\mathbb{Z}u+\mathbb{Z} \alpha u$ where $u=
\alpha-1$. Moreover
 $\partial{\mathcal R}_a = \bigcup_{v \in B} \mathcal{R}_a \cap (\mathcal{R}_a + v)$, where
  $B= \{\pm u, \pm  \alpha u, \pm(1+\alpha) u, \pm(\alpha-1) u \}$
and $ \mathcal{R}_a \cap (\mathcal{R}_a +  (1+ \alpha) u))=
\{-1\},\;
 \mathcal{R}_a \cap (\mathcal{R}_a +  ( \alpha-1) u)= \{-\alpha\}$.
\end{prop}
\begin{figure}[h!]
\centering \epsfig{file=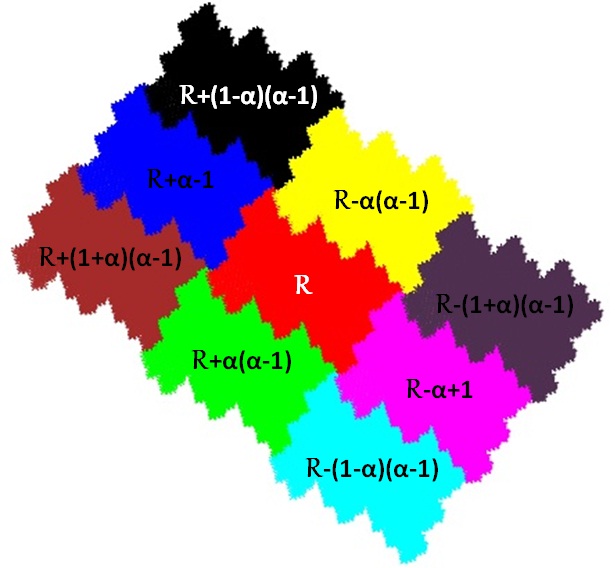,height=5cm,width=5.5cm}
\caption{ {\small Tiling induced by $\mathcal{R}_3$ }}
\end{figure}
\begin{obs}
 In this paper we will use the following relations:
\begin{equation}\label{rel}
\alpha^{n}-a\alpha^{n-1}+\alpha^{n-2}-\alpha^{n-3}=0,\;
\alpha^{n}=(a-1)\alpha^{n-1}+(a-1)\alpha^{n-2}+\alpha^{n-4},\forall\
n \in \mathbb{Z}.
\end{equation}
\end{obs}
\begin{lem}\label{ex1}
Let $z \in \mathcal{B}_{\alpha-1}$ then $z=(\alpha -1) +
\displaystyle\sum_{i=2}^{\infty}a_i\alpha^{i}$ and
$z=\displaystyle\sum_{i=2}^{\infty}b_i\alpha^{i}$ with $a_{2}=0,\
b_{2}=a-1$ and $b_{4}=0$.
\end{lem}
\begin{proof} Take $z\in\mathcal{B}_{\alpha-1}.$ Using relation (\ref{rel}) we have
$z=\displaystyle\alpha^{-3}+(a-1)\alpha^{-1}+(a-2)+
\sum_{i=2}^{\infty}a_i\alpha^{i}$ and
$z=\displaystyle\sum_{i=2}^{\infty}b_i\alpha^{i}$. Then the
admissible path, starting from $0$, in the automaton associated to
$z$ is
$$(1,0)(0,0)(a-1,0)(a-2,0)(0,0)(a_2,b_2)(a_3,b_3)....$$ Using the
automaton we have that $a_2=0$ , $b_2=a-1$ and then we can write
$z=\displaystyle\alpha^{-3}+(a-1)\alpha^{-1}+(a-2)+
\sum_{i=3}^{\infty}a_i\alpha^{i}=\alpha-1+\sum_{i=3}^{\infty}a_i\alpha^{i}$
and $z=(a-1)\alpha^2+\displaystyle\sum_{i=3}^{\infty}b_i\alpha^{i}$.
We also have that $b_4=0$.
\end{proof}

\section{Parametrization of $\partial R_a,\ a\geq 3$}

In this section  we give a complete description of
$\partial\mathcal{R}_a,\ a\geq 3$. By Proposition (\ref{prop35}) we
have that $\partial \mathcal{R}_a = \bigcup_{v \in B} \mathcal{R}_a
\cap (\mathcal{R}_a+v)$ where
$B=\{\pm(\alpha^{-3}+\alpha^{-1})=\pm(\alpha-1),\ \pm
(\alpha^{2}-\alpha),\ \pm (\alpha^{2}-1),\ \pm (\alpha-1)^2\}$.
Since $\mathcal{R}_a \cap (\mathcal{R}_a \pm v)$ is a point if $v=
\alpha^{2}-1$ or $v= (\alpha-1)^2$, we will study the others four
regions $ \mathcal{B}_v = \mathcal{R}_a \cap (\mathcal{R}_a + v)$
where $v \in \{\pm (\alpha-1),  \pm (\alpha^{2}-\alpha)\}$. For this
we will use the set $ \mathcal{B}_{\alpha-1} = \mathcal{R}_a \cap
(\mathcal{R}_a + \alpha-1)$ described by
$$\mathcal{B}_{\alpha-1}=\left\{z\in\mathbb{C},z=\alpha-1+\sum_{i=3}^{\infty}
a_i\alpha^{i}=(a-1)\alpha^2+\sum_{i=3}^{\infty}
b_i\alpha^{i}\right\}.
$$ In particular, we will prove the following results.
\begin{prop}\label{t1}
    Let $f_i, \ i=1,2,3$, be the functions defined by $f_1(z)=
    \alpha^{-1}-1+\alpha^{-1}z,\ f_2(z)=-(a-1)\alpha+\alpha^{-1}z$ and
    $f_3(z)=1-\alpha+z$. Then we have the following properties:
    \begin{enumerate}
        \item $\mathcal{B}_{\alpha^2-\alpha}=f_1(\mathcal{B}_{\alpha-1}),$
        \item $\mathcal{B}_{\alpha-\alpha^2}= f_2(\mathcal{B}_{\alpha-1}),$
        \item $\mathcal{B}_{1-\alpha}=f_3(\mathcal{B}_{\alpha-1}).$
        \item $\mathcal{B}_{\alpha-1}\cap f_1(\mathcal{B}_{\alpha-1}) = \mathcal{B}_{\alpha-1} \cap \mathcal{B}_{\alpha^2-\alpha} = \{-1\}.$
        \item $f_{1}(\mathcal{B}_{\alpha-1})\cap f_3(\mathcal{B}_{\alpha-1}) = \mathcal{B}_{\alpha^{2}-\alpha} \cap \mathcal{B}_{1-\alpha} = \{-\alpha\}.$
        \item $f_{2}(\mathcal{B}_{\alpha-1})\cap f_3(\mathcal{B}_{\alpha-1}) = \mathcal{B}_{\alpha-\alpha^{2}} \cap \mathcal{B}_{1-\alpha} = \{-\alpha^{2}\}.$
       \item $\mathcal{B}_{\alpha-1}\cap f_2(\mathcal{B}_{\alpha-1}) = \mathcal{B}_{\alpha-1} \cap \mathcal{B}_{\alpha-\alpha^{2}} = \{-(a-1)\alpha-\alpha^{-1}\}.$
       \item $\mathcal{B}_{\alpha-1}\cap f_3(\mathcal{B}_{\alpha-1}) = \mathcal{B}_{\alpha-1} \cap \mathcal{B}_{1-\alpha} = \emptyset.$
       \item $f_{1}(\mathcal{B}_{\alpha-1})\cap f_2(\mathcal{B}_{\alpha-1}) = \mathcal{B}_{\alpha^{2}-\alpha} \cap \mathcal{B}_{\alpha-\alpha^{2}} = \emptyset.$

  \end{enumerate}
\end{prop}
\begin{prop}\label{t2} Let $g_i,\ i=0,1,...,2(a-1),$ be the functions defined by $g_{2k+1}(z)=-1-k\alpha^3+\alpha^3z$ for
    $k=0,...,a-2$, and $g_{2k}(z)=\alpha-1+(a-1-k)\alpha^3+\alpha^2z$
    for $k=0,...,a-1$. Then
    $$\displaystyle
    \mathcal{B}_{\alpha-1}=\bigcup_{i=0}^{2(a-1)}g_i(X_i),$$ where
    $X_i=\mathcal{B}_{\alpha-1}$ if $i$ is an odd number or $i=2(a-1)$
    and $X_i=\mathcal{B}'_{\alpha-1}=\{z\in \mathcal{B}_{\alpha-1}; a_3\neq a-1\}$
    if $i$ is an even number.
\end{prop}

\begin{obs}

Using  Proposition \ref{t2} we will construct an explicit continuous
and bijective application from $[0,1]$ to $\mathcal{R}_{\alpha-1}$.
Using this fact and Proposition \ref{t1} we obtain an explicit
homeomorphism between the circle and the boundary of ${\mathcal
R}_a$.
\end{obs}
$\mathbf{Proof\ of\  proposition\  \ref{t1}:}$ According to
(\ref{rel}) we have $\alpha^2-\alpha=\alpha^{-2}+(a-1)+(a-2)\alpha$.
\begin{enumerate}
    \item Take $z\in \mathcal{B}_{\alpha-1}$. According to Lemma \ref{ex1} $ \displaystyle z=\alpha-1+\sum_{i=3}^{\infty}
    a_i\alpha^{i}=(a-1)\alpha^2+\sum_{i=3}^{\infty}
    b_i\alpha^{i}$. Then \\
    $\begin{array}{l}f_1(z)=\displaystyle\alpha^{-1}-1+\alpha^{-1}((a-1)\alpha^2+\sum_{i=3}
    b_i\alpha^{i}) =\alpha^{-2}+(a-1)+(a-2)\alpha+\sum_{i=3}
    b_i\alpha^{i-1}\in
    \mathcal{R}_a+\alpha^2-\alpha.\end{array}$\\ We also have\\
    $\begin{array}{l}\displaystyle
    f_1(z)=\alpha^{-1}-1+\alpha^{-1}(\alpha-1+\sum_{i=3}
    a_i\alpha^{i})
    =\sum_{i=3} a_i\alpha^{i-1}\in \mathcal{R}_a.\end{array}$\\ Therefore
    $f_1(\mathcal{B}_{\alpha-1})\subseteq
    \mathcal{B}_{\alpha^2-\alpha}.$ \\
    Take $z\in \mathcal{B}_{\alpha^2-\alpha},\ \displaystyle
    z=\alpha^{-2}+(a-1)+(a-2)\alpha+\sum_{i=2}
    a_i\alpha^{i}=\sum_{i=2} b_i\alpha^{i}.$ Then\\
    $\begin{array}{l}\displaystyle
    f^{-1}_1(z)=\alpha-1+\alpha(\alpha^{-2}+(a-1)+(a-2)\alpha+\sum_{i=2}
    a_i\alpha^{i})=(a-1)\alpha^2+\sum_{i=2} a_i\alpha^{i+1}\in \mathcal{R}_a.\end{array}$\\ We also have\\
    $\begin{array}{l}\displaystyle
    f^{-1}_1(z)=\alpha-1+\alpha(\sum_{i=2}
    b_i\alpha^{i})=\alpha-1+\sum_{i=2} b_i\alpha^{i+1}\in
    \mathcal{R}_a+\alpha-1.\end{array}$\\
    Therefore $f^{-1}_1(\mathcal{B}_{\alpha^2-\alpha})\subseteq
    \mathcal{B}_{\alpha-1}$  and then
    $$f_1(\mathcal{B}_{\alpha-1})= \mathcal{B}_{\alpha^2-\alpha}.$$
    \item Take $z\in\mathcal{B}_{\alpha-\alpha^2}$. Then
    $z+\alpha^2-\alpha$ belongs to $\mathcal{B}_{\alpha^2-\alpha}$ and
    according to what was done before there exists
    $w\in\mathcal{B}_{\alpha-1}$ such that $z+\alpha^2-\alpha=g_1(w)$.
    Then\\
    $z+\alpha^2-\alpha=\alpha^{-1}-1+\alpha^{-1}(w)\Rightarrow
    z=\alpha^{-1}-1+\alpha-\alpha^2+\alpha^{-1}(w)=-(a-1)\alpha+\alpha^{-1}(w).$
    Therefore
    $$f_2(\mathcal{B}_{\alpha-1})= \mathcal{B}_{\alpha-\alpha^2}.$$ We
    also know that
    $f_1^{-1}(z+\alpha^2-\alpha)=\alpha-1+\alpha(z+\alpha^2-\alpha)=(a-1)\alpha^2+\alpha
    z=f_2^{-1}(z)\in \mathcal{B}_{\alpha-1}$ and then
    $f_2^{-1}(\mathcal{B}_{\alpha-\alpha^2})\subseteq
    \mathcal{B}_{\alpha-1} $. Therefore $f_2(\mathcal{B}_{\alpha-1})=
    \mathcal{B}_{\alpha-\alpha^2}.$
    \item This item can be done by the same manner of item $(2)$.

    \item Take $z \in \mathcal{B}_{\alpha-1} \cap \mathcal{B}_{\alpha^2-\alpha} = \mathcal{R} \cap (\mathcal{R} + \alpha -1) \cap (\mathcal{R} + \alpha^{2} - \alpha).$ \\Then $z-\alpha + 1 \in \mathcal{R}  \cap (\mathcal{R} - \alpha + 1) \cap (\mathcal{R} + (\alpha - 1)^{2}) \subseteq \mathcal{R} \cap (\mathcal{R} + (\alpha-1)^{2}) = \{-1\}$. \\Therefore, $z-\alpha+1=-\alpha,$ and $z=-1$.




    \item [(8)]Take $z \in \mathcal{B}_{\alpha-1} \cap f_{3}(\mathcal{B}_{\alpha-1})$. Then there is $z_{1} \in \mathcal{B}_{\alpha -1}$ such that $z = 1-\alpha+z_{1}$. Then
    \[
    \alpha + z = 1 + z_{1}.
    \]
    Since $z, \ z_{1} \in \mathcal{B}_{\alpha -1}$ we know that $z=(a-1)\alpha^{2} + \sum_{i=3}^{\infty}a_{i}\alpha^{i}$ and $z_{1} = (a-1)\alpha^{2} + \sum_{i=3}^{\infty}b_{i}\alpha^{i}.$

    Then the equality above becomes
    \[
    1+(a-1)\alpha^{2} + \sum_{i=3}^{\infty}a_{i}\alpha^{i}= \alpha +  (a-1)\alpha^{2} + \sum_{i=3}^{\infty}b_{i}\alpha^{i}.
    \]

  So, we conclude that $(1,0)(0,1)(a-1,a-1) \ldots$ is an admissible path on the automaton $\mathcal{A}$ starting from $0$. But there is no such path on the automaton and then
  \[
  \mathcal{B}_{\alpha-1}\cap f_3(\mathcal{B}_{\alpha-1}) = \mathcal{B}_{\alpha-1} \cap \mathcal{B}_{1-\alpha} = \emptyset.
  \]



Following the ideas of  items $(4)$ and $(8)$ we can prove $(5), (6), (7) $ and $(9)$. \\For more details see \cite{JT}

    \hfill$\blacksquare$

\begin{figure}[h]
  \centering \epsfig{file=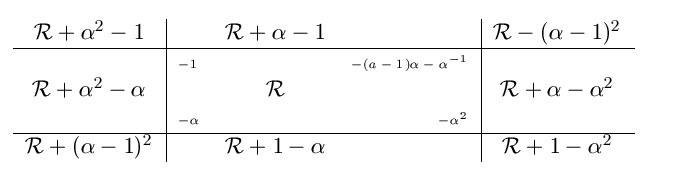,height=4cm,width=12cm}
   \caption{
      {\small Boundary of $\mathcal{R}_{a}$ }}\label{F4}
\end{figure}

\end{enumerate}

$\mathbf{Proof \ of\ Proposition\ \ref{t2}:}$ Let $z$ be an element
of $\mathcal{B}_{\alpha-1}$. Using the automaton $\mathcal{A}$ we
can write
$z=\alpha-1+\displaystyle\sum_{i=3}a_i\alpha^{i}=(a-1)\alpha^2+\displaystyle\sum_{i=3}b_i\alpha^{i}$
where $(a_3,b_3)=(t,t), \ t=0,1,...,a-1$ or $(a_3,b_3)=(t,t-1), \
t=1,...,a-1$. Let $\mathcal{B}_{\alpha-1}^{1,t},\
\mathcal{B}_{\alpha-1}^{2,t}$ be the
following sets:\\
$\mathcal{B}_{\alpha-1}^{1,t}=\{z\in \mathcal{B}_{\alpha-1};
(a_3,b_3)=(t,t),\ t=0,1...,a-1\}$, \\ \
\\$\mathcal{B}_{\alpha-1}^{2,t}=\{z\in \mathcal{B}_{\alpha-1};
(a_3,b_3)=(t,t-1),\ t=1,2...,a-1\}$. Since
$$\mathcal{B}_{\alpha-1}=\left[\bigcup_{t=0}^{a-1}\mathcal{B}_{\alpha-1}^{1,t}\right]\bigcup\left[\bigcup_{t=1}^{a-1} \mathcal{B}_{\alpha-1}^{2,t}\right],$$
in order to prove this theorem we need to show that
$g_{2k+1}(\mathcal{B}_{\alpha-1})=\mathcal{B}_{\alpha-1}^{2,a-1-k},
k=0,...,a-2$,\\
$g_{2(a-1)}(\mathcal{B}_{\alpha-1})=\mathcal{B}_{\alpha-1}^{1,0}$
and
$g_{2k}(\mathcal{B}_{\alpha-1}^{'})=\mathcal{B}_{\alpha-1}^{1,a-1-k}, k=0,...,a-2$.\\
1)- Indeed since $z\in\mathcal{B}_{\alpha-1}$ then: \\
$g_{2k+1}(z)=-1-k\alpha^3+\alpha^3(\alpha-1+\displaystyle
\sum_{i=3}a_i\alpha^{i})=(a-1)\alpha^2+(a-2-k)\alpha^3+\sum_{i=3}a_i\alpha^{i+3}$
and \\
$g_{2k+1}(z)=-1-k\alpha^3+\alpha^3((a-1)\alpha^2+\displaystyle
\sum_{i=3}b_i\alpha^{i})=\alpha-1+(a-1-k)\alpha^3+(a-1)\alpha^4+(a-2)\alpha^5+\displaystyle\sum_{i=3}b_i\alpha^{i+3}$,
that is $g_{2k+1}(\mathcal{B}_{\alpha-1})\subseteq
\mathcal{B}_{\alpha-1}^{2,a-1-k}$. On the other hand if we take
$w\in \mathcal{B}_{\alpha-1}^{2,a-1-k}$,
$w=\alpha-1+(a-1-k)\alpha^3+(a-1)\alpha^4+(a-2)\alpha^5+\displaystyle\sum_{i=6}u_i\alpha^{i}=(a-1)\alpha^2+(a-2-k)\alpha^3+\displaystyle\sum_{i=6}v_i\alpha^{i}$
then
$z=\alpha-1+\displaystyle\sum_{i=6}v_i\alpha^{i-3}=(a-1)\alpha^2+\displaystyle\sum_{i=6}u_i\alpha^{i-3}$
is an element of $\mathcal{B}_{\alpha-1}$ such that $g_{2k+1}(z)=w$.
Therefore we conclude that
$$g_{2k+1}(\mathcal{B}_{\alpha-1})=\mathcal{B}_{\alpha-1}^{2,a-1-k},
k=0,...,a-2.$$
$g_{2(a-1)}(z)=\alpha-1+\alpha^2(\alpha-1+\displaystyle
\sum_{i=3}a_i\alpha^{i})=(a-1)\alpha^2+\sum_{i=3}a_i\alpha^{i+2}$
and \\
$g_{2(a-1)}(z)=\alpha-1+\alpha^2((a-1)\alpha^2+\displaystyle
\sum_{i=3}b_i\alpha^{i})=\alpha-1+(a-1)\alpha^4+\sum_{i=3}b_i\alpha^{i+2}$,
that is $g_{2(a-1)}(\mathcal{B}_{\alpha-1})\subseteq
\mathcal{B}_{\alpha-1}^{1,0}$. On the other hand if we take $w\in
R_{\alpha-1}^{1,0}$,
$w=\alpha-1+(a-1)\alpha^4+\displaystyle\sum_{i=5}u_i\alpha^{i}=(a-1)\alpha^2+\displaystyle\sum_{i=5}v_i\alpha^{i}$
then
$z=\alpha-1+\displaystyle\sum_{i=5}v_i\alpha^{i-2}=(a-1)\alpha^2+\displaystyle\sum_{i=5}u_i\alpha^{i-2}$
is an element of $\mathcal{B}_{\alpha-1}$ such that
$g_{2(a-1)}(z)=w$. Therefore
$$g_{2(a-1)}(\mathcal{B}_{\alpha-1})=\mathcal{B}_{\alpha-1}^{1,0}.$$
2)- Let $z\in \mathcal{B}'_{\alpha-1}$ given by
$z=\alpha-1+\displaystyle\sum_{i=3}a_i\alpha^{i}=(a-1)\alpha^2+\displaystyle\sum_{i=3}b_i\alpha^{i}$.
Since $a_3\neq a-1$ then using the automaton we have $b_3\neq a-1$.
Then\\
$g_{2k}(z)=\alpha-1+(a-1-k)\alpha^3+\alpha^2(\alpha-1+\displaystyle\sum_{i=3}a_i\alpha^{i})=(a-1)\alpha^2+(a-1-k)\alpha^3
+\displaystyle\sum_{i=3}a_i\alpha^{i+2}$ and\\
$g_{2k}(z)=\alpha-1+(a-1-k)\alpha^3+\alpha^2((a-1)\alpha^2+\displaystyle\sum_{i=3}b_i\alpha^{i})=\alpha-1+(a-1-k)\alpha^3+(a-1)\alpha^4
+\displaystyle\sum_{i=3}b_i\alpha^{i+2}.$ So we have
$g_{2k}(\mathcal{B}'_{\alpha-1})\subseteq
\mathcal{B}_{\alpha-1}^{1,a-1-k}$. On the other hand if we take
$w\in \mathcal{B}_{\alpha-1}^{1,a-1-k}$,
$w=\alpha-1+(a-1-k)\alpha^3+(a-1)\alpha^4+\displaystyle\sum_{i=5}u_i\alpha^{i}=(a-1)\alpha^2+(a-1-k)\alpha^3+\displaystyle\sum_{i=5}v_i\alpha^{i}$
then we have $u_5,v_5\neq a-1$ (again use the automaton) and
$z=\alpha-1+\displaystyle\sum_{i=5}v_i\alpha^{i-2}=(a-1)\alpha^2+\displaystyle\sum_{i=5}u_i\alpha^{i-2}$
is an element of $\mathcal{B}'_{\alpha-1}$ such that $g_{2k}(z)=w$.
Therefore
$$g_{2k}(\mathcal{B}_{\alpha-1}^{'})=\mathcal{B}_{\alpha-1}^{1,a-1-k},
k=0,...,a-2.$$\hfill$\blacksquare$

Using the previous notation and taking $u=-1,
v=-(a-1)\alpha-\alpha^{-1}, w=-1-\alpha^3$, we have the following
lemmas.
\begin{lem}\label{l1}
    \
    \begin{enumerate}
        \item Take $k_1\leq k_2$. Then $g_{2k_1}(\mathcal{B}'_{\alpha-1})\bigcap g_{2k_2+1}(\mathcal{B}_{\alpha-1})=\left\{\begin{array}{l}\emptyset,\mbox{ if } k_2>k_1\\ -1-\alpha^2-k\alpha^3-(a-1)\alpha^4, \mbox{ if } k_2=k_1.\end{array}\right.$\\ Therefore $-1-\alpha^2-k\alpha^3-(a-1)\alpha^4=g_{2k}(w)=g_{2k+1}(v)$.

        \item Take $k_1< k_2$.Then $g_{2k_1+1}(\mathcal{B}_{\alpha-1})\bigcap g_{2k_2}(\mathcal{B}'_{\alpha-1})=\left\{\begin{array}{l}\emptyset,\mbox{ if } k_2>k_1+1\\ -1-(k_1+1)\alpha^3, \mbox{ if } k_2=k_1+1.\end{array}\right.$\\ Therefore $-1-(k_1+1)\alpha^3=g_{2k+1}(u)=g_{2(k+1)}(v)$.

        \item $g_{2k_1}(\mathcal{B}'_{\alpha-1})\bigcap g_{2k_2}(\mathcal{B}'_{\alpha-1})=\emptyset, \mbox{ if } k_1\neq k_2$.
        \item $g_{2k_1+1}(\mathcal{B}_{\alpha-1})\bigcap g_{2k_2+1}(\mathcal{B}_{\alpha-1})=\emptyset, \mbox{ if } k_1\neq k_2$.
        \item $\displaystyle\lim_{n\longrightarrow \infty}(g_0\circ g_{2(a-1)})^n(z)=u=-1,\ \forall z\in \mathcal{B}_{\alpha-1}.$
        \item $\displaystyle\lim_{n\longrightarrow \infty}(g_{2(a-1)}\circ g_0)^n(z)=v=-(a-1)\alpha-\alpha^{-1},\ \forall z\in\mathcal{B}'_{\alpha-1}.$
    \end{enumerate}
\end{lem}
\begin{dem}
    $1)-$ Take $z\in g_{2k_1}(\mathcal{B'}_{\alpha-1})\bigcap g_{2k_2+1}(\mathcal{B}_{\alpha-1})$, $k_1\leq k_2$. \\Then $z=g_{2k_1}(z_1)=g_{2k_2+1}(z_2)$, $ z_1\in \mathcal{B'}_{\alpha-1}, z_2\in \mathcal{B}_{\alpha-1}$ and $$\alpha -1+(a-1-k_1)\alpha^3+\alpha^2 z_1=-1-k_2\alpha^3+\alpha^3z_2.$$ If we suppose $k_2=k_1+k$ we have $$\alpha+(a-1)\alpha^3+\alpha^2z_1=-k\alpha^3+\alpha^3z_2,$$ and multiplying by $\alpha^{-3}$
    $$\alpha^{-2}+(a-1)+\alpha^{-1}z_1=-k+z_2.$$ Since $ z_1\in \mathcal{B'}_{\alpha-1}, z_2\in \mathcal{B}_{\alpha-1}$ we know that $z_1=\alpha-1+\sum_{i=3}^{\infty}a_{i}\alpha^i, a_3\neq a-1$ and $z_2=\alpha-1+\sum_{i=3}^{\infty}b_{i}\alpha^i.$ Then the equality above becomes\\
    $$ \alpha^{-2}+a-\alpha^{-1}+\sum_{i=3}^{\infty}a_{i}\alpha^{i-1}=-k+\alpha-1 +\sum_{i=3}^{\infty}b_{i}\alpha^i, $$ and since $\alpha^{-2}+a-1=\alpha$ then $$(k+1) +\sum_{i=3}^{\infty}a_{i}\alpha^{i-1}=\sum_{i=3}^{\infty}b_{i}\alpha^i.$$
    So we conclude that $(k+1,0)(0,0)(a_3,0)(a_4,b_3)(a_5,b_4)...$ is an admissible path on the automaton $\mathcal{A}$ starting from $0$. Using the automaton, since $a_3\neq a-1$, we see that the only possibility is $$(1,0)(0,0)(a-2,0)(a-1,0)(0,a-1)(0,a-1)(a-1,0)(a-1,0)...$$
    Then $k=0$,  $$z_1=\alpha-1+(a-2)\alpha^3+(a-1)\alpha^4+\sum_{i=2}^{\infty}[(a-1)\alpha^{4i-1}+(a-1)\alpha^{4i}]=-1-\alpha^3,$$and $$z_2=\alpha-1+\sum_{i=1}^{\infty}[(a-1)\alpha^{4i}+(a-1)\alpha^{4i+1}]=-(a-1)\alpha-\alpha^{-1}.$$
    $2)-$   Take $z\in g_{2k_1+1}(\mathcal{B}_{\alpha-1})\bigcap g_{2k_2}(\mathcal{B'}_{\alpha-1})$, $k_1< k_2$. \\Then $z=g_{2k_1+1}(z_1)=g_{2k_2}(z_2)$, $ z_1\in \mathcal{B}_{\alpha-1}, z_2\in \mathcal{B'}_{\alpha-1}$ and $$-1-k_1\alpha^3+\alpha^3z_1=\alpha -1+(a-1-k_2)\alpha^3+\alpha^2 z_2.$$ If we suppose $k_2=k_1+k$ we have $$\alpha^3z_1=\alpha+(a-1-k)\alpha^3+\alpha^2 z_2,$$ and multiplying by $\alpha^{-3}$
    $$z_1=\alpha^{-2}+(a-1-k)+\alpha^{-1}z_2.$$ Since $ z_1\in \mathcal{B}_{\alpha-1}, z_2\in \mathcal{B'}_{\alpha-1}$ we know that $z_1=(a-1)\alpha^{2}+\sum_{i=3}^{\infty}a_{i}\alpha^i$ and $z_2=(a-1)\alpha^{2}+\sum_{i=3}^{\infty}b_{i}\alpha^i.$ Then the equality above becomes\\
    $$ (a-1)\alpha^{2}+\sum_{i=3}^{\infty}a_{i}\alpha^i=\alpha^{-2} +(a-1-k)+(a-1)\alpha+\sum_{i=3}^{\infty}b_{i}\alpha^{i-1}.$$
    So we conclude that $(1,0)(0,0)(a-1-k,0)(a-1,0)(b_3,a-1)(b_4,a_3)(b_5,a_4)...$ is an admissible path on the automaton $\mathcal{A}$ starting from $0$. Using the automaton, we see that the only possibility is $$(1,0)(0,0)(a-2,0)(a-1,0)(0,a-1)(0,a-1)(a-1,0)(a-1,0)...$$
    Then $k=1$,  $$z_2=(a-1)\alpha^2+\sum_{i=1}^{\infty}[(a-1)\alpha^{4i+1}+(a-1)\alpha^{4i+2}]=-(a-1)\alpha-\alpha^{-1},$$and $$z_1=\sum_{i=1}^{\infty}[(a-1)\alpha^{4i-2}+(a-1)\alpha^{4i-1}]=-1.$$
    $3)-$Take $z\in g_{2k_1}(\mathcal{B'}_{\alpha-1})\bigcap g_{2k_2}(\mathcal{B'}_{\alpha-1})$, $k_1< k_2$.\\ Then $z=g_{2k_1}(z_1)=g_{2k_2}(z_2)$, $ z_1\in \mathcal{B'}_{\alpha-1}, z_2\in \mathcal{B'}_{\alpha-1}$ and $$\alpha-1+(a-1-k_1)\alpha^3+\alpha^2z_1=\alpha -1+(a-1-k_2)\alpha^3+\alpha^2 z_2.$$ If we suppose $k_2=k_1+k$ we have $$k\alpha^3+\alpha^2z_1=\alpha^2 z_2,$$ and multiplying by $\alpha^{-2}$
    $$k\alpha+z_1=z_2.$$ Since $ z_1\in \mathcal{B'}_{\alpha-1}, z_2\in \mathcal{B'}_{\alpha-1}$ we know that $z_1=(a-1)\alpha^{2}+\sum_{i=3}^{\infty}a_{i}\alpha^i$, $z_2=(a-1)\alpha^{2}+\sum_{i=3}^{\infty}b_{i}\alpha^i,$ and, by Example \ref{ex1}, $a_4=b_4=0.$ Then the equality above becomes\\
    $$k\alpha+ (a-1)\alpha^{2}+\sum_{i=3}^{\infty}a_{i}\alpha^i=(a-1)\alpha^2+\sum_{i=3}^{\infty}b_{i}\alpha^{i}.$$
    So we conclude that if the intersection is not empty, $(k,0)(a-1,a-1)(a_3,b_3)(0,0)(a_5,b_5)...$ is an admissible path on the automaton $\mathcal{A}$ starting from $0$. But there is no such path on the automaton and then $$g_{2k_1}(\mathcal{B'}_{\alpha-1})\bigcap g_{2k_2}(\mathcal{B'}_{\alpha-1})=\emptyset.$$
    Using the same ideas we can prove that $g_{2k_1+1}(\mathcal{B}_{\alpha-1})\bigcap g_{2k_2+1}(\mathcal{B}_{\alpha-1})=\emptyset$.\\
    5)- Using induction we can prove that
    \[(g_0 \circ g_{r-1})^n(z)=\sum_{i=1}^{n}\left[(a-1)\alpha^{4i-2}+(a-1)\alpha^{4i-1}\right]+\alpha^{4n}z.\]
    Then $\displaystyle\lim_{n\longrightarrow \infty}(g_0\circ
    g_{r-1)})^n(z)=\sum_{i=1}^{\infty}\left[(a-1)\alpha^{4i-2}+(a-1)\alpha^{4i-1}\right]=-1.$\\
    Indeed for $n=1$ and $z\in \mathcal{B}_{\alpha-1}$ by proposition \ref{t2} we
    have $g_{r-1}(z)=\alpha-1+\alpha^2z\in \mathcal{B}_{\alpha-1}^{1,0}$.
    By definition we have\\
    $\begin{array}{ll}g_0(g_{r-1}(z))&=\alpha-1+(a-1)\alpha^3+\alpha^2(\alpha-1+\alpha^2z)=\alpha-1+(a-1)\alpha^3+\alpha^3-\alpha^2+\alpha^4z=\\
    \ &=(a-1)\alpha^2+(a-1)\alpha^3+\alpha^4z.\end{array}$

    Suppose the formula is true for $k \geq 1$, that is
    \[(g_0\circ g_{r-1})^k(z)=\sum_{i=1}^{k}\left[(a-1)\alpha^{4i-2}+(a-1)\alpha^{4i-1}\right]+\alpha^{4k}z.\
    \ \ (*)\]
    We have to prove the formula for $n=k+1$.
    Since $(g_0\circ g_{r-1})^{k+1}(z)=(g_0\circ g_{r-1})\circ(g_0\circ g_{r-1})^{k}(z)$ using $(*)$ we have to prove that\\
    $\displaystyle g_0(g_{r-1}(\sum_{i=1}^{k}\left[(a-1)\alpha^{4i-2}+(a-1)\alpha^{4i-1}\right]+\alpha^{4k}z))= \sum_{i=1}^{k+1}\left[(a-1)\alpha^{4i-2}+(a-1)\alpha^{4i-1}\right]+\alpha^{4(k+1)}z.$
    Indeed
    \\${\small \begin{array}{l}g_{r-1}(\displaystyle\sum_{i=1}^{k}\left[(a-1)\alpha^{4i-2}+(a-1)\alpha^{4i-1}\right]+\alpha^{4k}z)=\alpha-1+\alpha^2(\sum_{i=1}^{k}\left[(a-1)\alpha^{4i-2}+(a-1)\alpha^{4i-1}\right]+\alpha^{4k}z)=\\
        \alpha-1+\displaystyle\sum_{i=1}^{k}\left[(a-1)\alpha^{4i}+(a-1)\alpha^{4i+1}\right]+\alpha^{4k+2}z,\end{array}}$\\
    and\\
    ${\small
        \begin{array}{l}g_0(g_{r-1}(\displaystyle\sum_{i=1}^{k}\left[(a-1)\alpha^{4i-2}+(a-1)\alpha^{4i-1}\right]+\alpha^{4k}z))=
        g_0(\alpha-1+\displaystyle\sum_{i=1}^{k}\left[(a-1)\alpha^{4i}+(a-1)\alpha^{4i+1}\right]+\alpha^{4k+2}z)=\\
        =\alpha-1+(a-1)\alpha^3+\alpha^2(\alpha-1+\displaystyle\sum_{i=1}^{k}\left[(a-1)\alpha^{4i}+(a-1)\alpha^{4i+1}\right]+\alpha^{4k+2}z)=\\
        =
        (a-1)\alpha^2+(a-1)\alpha^3+\displaystyle\sum_{i=1}^{k}\left[(a-1)\alpha^{4i+2}+(a-1)\alpha^{4i+3}\right]+\alpha^{4(k+1)}z
        =\displaystyle\sum_{i=1}^{k+1}\left[(a-1)\alpha^{4i-2}+(a-1)\alpha^{4i-1}\right]+\alpha^{4(k+1)}z.\end{array}}$\\
    \ \\
    $6)-$ Using induction we can prove that $$(g_{r-1}\circ
    g_0)^n(z)=(a-1)\alpha^2+\displaystyle\sum_{i=1}^{n-1}\left[(a-1)\alpha^{4i+1}+(a-1)\alpha^{4i+2}\right]+(a-1)\alpha^{4n+1}+\alpha^{4n}z,$$ and then
    $\displaystyle\lim_{n\longrightarrow \infty}(g_{r-1}\circ
    g_0)^n(z)=(a-1)\alpha^2+\displaystyle\sum_{i=1}^{\infty}\left[(a-1)\alpha^{4i+1}+(a-1)\alpha^{4i+2}\right]=-(a-1)\alpha-\alpha^{-1}.$

    \hfill$\blacksquare$

\end{dem}
\vspace{0.3cm}

%
%
\begin{prop}\label{l3}
Let $t \in [0,1],a\geq 3, r=2a-1$. Then there exists an unusual
expansion $(a_i)_{i\geq1}\in\{0,1,...,r-1\}^{\mathbb{N}}$,
$n_i,m_i\in\mathbb{N}$ such that we can write
$$t=\displaystyle\frac{a_{1}}{r}+\sum_{{i=2} \atop
{m_{i}+n_{i}=k}}^{\infty}
\displaystyle\frac{a_{i}}{r^{n_{i}}(r-2)^{m_{i}}},$$ where the
digits $a_i$ and the numbers $n_i,m_i$ satisfy the following
properties:
\begin{enumerate}

\item if
$a_1\in\{1,3,5,...,r-2,r-1\}$ then $a_2\in\{0,1,...,r-1\},\ n_2=2,\
m_2=0$;
\item if
$a_1\in\{0,2,4,...,r-3\}$ then $a_2\in\{0,1,...,r-3\},\ n_2=1,\ m_2=1$;\\
\ and for $i\geq 3$ we have:

\item $a_{i} \in \{0,1,2...,r-1\},\ m_i=m_{i-1},\ n_i=n_{i-1}+1$ if one of the following conditions are satisfied
\begin{enumerate}
\item $a_{i-1}=0$ and $i-1$ even;
\item $a_{i-1}=r-1$ and $i-1$ odd;
\item $a_{i-1}=2n-1,\ n=1,...,a-1$;
\item $a_{i-1}=r-3$, $i-1$ odd, $n_{i-1}=n_{i-2}$ and $m_{i-1}=m_{i-2}+1$
\end{enumerate}

\item $a_{i} \in \{0,1,2...,r-3\},\ m_i=m_{i-1}+1,\ n_i=n_{i-1}$ if one of the following conditions are satisfied
\begin{enumerate}
\item $a_{i-1}=0$ and $i-1$ odd;
\item $a_{i-1}=r-1$ and $i-1$ even;
\item $a_{i-1}=2n,\ n=1,...,a-3$ or $a_{i-1}=r-3$, $i-1$ even or odd, $n_{i-1}=n_{i-2}+1$ and
$m_{i-1}=m_{i-2}$;
\item $a_{i-1}=r-3$, $i-1$ even, $n_{i-1}=n_{i-2}$ and $m_{i-1}=m_{i-2}+1$
\end{enumerate}
\end{enumerate}
\end{prop}
\dem We can write $$1=\frac{r-1}{r}+\sum_{i=0}^{\infty}
\frac{r-1}{r^{2+i}(r-2)^{i}}+\frac{r-3}{r^{2+i}(r-2)^{1+i}}.
$$Given $t\in[0,1)$, we can prove by induction that for each $k\geq
1$ we can write
$$t=\displaystyle\frac{a_{1}}{r} +\sum_{{i=2} \atop
{m_{i}+n_{i}=i}}^{k}
\displaystyle\frac{a_{k}}{r^{n_{i}}(r-2)^{m_{i}}}+c_k,\ 0\leq
c_k<\frac{1}{r^{n_{k}}(r-2)^{m_{k}}}.$$Indeed if $t\in[0,1)$, then
there exist $a_1\in\{0,1,...,r-1\}$ such that $rt=a_1+t_1$, with
$0\leq t_1<1$. Then
$t=\frac{a_1}{r}+\frac{t_1}{r}=\frac{a_1}{r}+c_1$ with
$c_1<\frac{1}{r}$.
\\ If $a_1\in\{1,3,...,r-2,r-1\}$ then there exist
$a_2\in\{0,1,...,r-1\}$ such that $rc_1=\frac{a_2}{r}+t_2$, with
$0\leq t_2<\frac{1}{r}$ and then
$$t=\frac{a_1}{r}+c_1=\frac{a_1}{r}+\frac{a_2}{r^2}+\frac{t_2}{r}=\frac{a_1}{r}+\frac{a_2}{r^2}+c_2,\ 0<c_2<\frac{1}{r^2}.$$
If $a_1\in\{0,2,...,r-3\}$ then there exist $a_2\in\{0,1,...,r-3\}$
such that $rc_1=\frac{a_2}{r-2}+t_2$, with $0\leq t_2<\frac{1}{r-2}$
and then
$$t=\frac{a_1}{r}+c_1=\frac{a_1}{r}+\frac{a_2}{r(r-2)}+\frac{t_2}{r}=\frac{a_1}{r}+\frac{a_2}{r(r-2)}+c_2,\ 0<c_2<\frac{1}{r(r-2)}.$$
Then the result is true for $k=1,2$. Suppose that it is true for
$2<k$, that is $$t=\displaystyle\frac{a_{1}}{r}+\sum_{{i=2} \atop
{m_{i}+n_{i}=i}}^{k}
\displaystyle\frac{a_{k}}{r^{n_{i}}(r-2)^{m_{i}}}+c_k,
$$ where $0\leq c_k<\frac{1}{r^{n_{k}}(r-2)^{m_{k}}}$.\\ If $a_k$
and $k$ satisfy condition $(3a)$ or $(3b)$ or $(3c)$ or $(3d)$ then
there exist
$a_{k+1}\in\{0,1,...,r-1\}$ such that\\
$r^{n_{k}}(r-2)^{m_{k}}c_k=\frac{a_{k+1}}{r}+t_{k+1},\ 0\leq
t_{k+1}<\frac{1}{r}$ and then $$t=\displaystyle\frac{a_{1}}{r}
+\sum_{{i=2} \atop {m_{i}+n_{i}=i}}^{k}
\displaystyle\frac{a_{k}}{r^{n_{i}}(r-2)^{m_{i}}}+c_k=\displaystyle\frac{a_{1}}{r}
+\sum_{{i=2} \atop {m_{i}+n_{i}=i}}^{k}
\displaystyle\frac{a_{k}}{r^{n_{i}}(r-2)^{m_{i}}}+
\frac{a_{k+1}}{r^{n_k+1}(r-2)^{m_{k}}}+c_{k+1},$$ where
$0\leq c_{k+1}=\displaystyle\frac{t_{k+1}}{r^{n_{k}}(r-2)^{m_{k}}}<\frac{1}{r^{n_{k}+1}(r-2)^{m_{k}}}.$ Then the result is true for $k+1$.\\
If $a_k$ and $k$ satisfy condition $(4a)$ or $(4b)$ or $(4c)$ or
$(4d)$ then there exist
$a_{k+1}\in\{0,1,...,r-3\}$ such that\\
$r^{n_{k}}(r-2)^{m_{k}}c_k=\frac{a_{k+1}}{r-2}+t_{k+1},\ 0\leq
t_{k+1}<\frac{1}{r-2}$ and then
$$t=\displaystyle\frac{a_{1}}{r}+\sum_{{i=2} \atop
{m_{i}+n_{i}=i}}^{k}
\displaystyle\frac{a_{k}}{r^{n_{i}}(r-2)^{m_{i}}}+c_k=\displaystyle\frac{a_{1}}{r}
+\sum_{{i=2} \atop {m_{i}+n_{i}=i}}^{k}
\displaystyle\frac{a_{k}}{r^{n_{i}}(r-2)^{m_{i}}}+
\frac{a_{k+1}}{r^{n_k}(r-2)^{m_{k}+1}}+c_{k+1},$$ where $0\leq
c_{k+1}=\displaystyle\frac{t_{k+1}}{r^{n_{k}}(r-2)^{m_{k}}}<\frac{1}{r^{n_{k}}(r-2)^{m_{k}+1}}.$
Then the result is true for $k+1$. Therefore the result is true for
every $k\geq1$. $\hfill$$\blacksquare$
\begin{obs}
Let $t\in[0,1]$ written as
$t=\displaystyle\frac{a_{1}}{r}+\sum_{{i=2} \atop
{m_{i}+n_{i}=i}}^{k}
\displaystyle\frac{a_{k}}{r^{n_{i}}(r-2)^{m_{i}}}$. In order to
simplify the demonstration of some of the results in this paper, a
simpler notation will be used, i. e., $t$ wil be represented as
$$t=\frac{a_1}{r}+\sum_{{i=2} \atop
{m_{i}+n_{i}=i}}^{\infty}\frac{a_i}{(n_i,m_i)}. $$
\end{obs}
    \begin{prop}\label{p37}
        Let $t,t'$ be elements in $[0,1]$,
        $t=\frac{a_1}{r}+\sum_{i=2}^{\infty}\frac{a_i}{(n_i,m_i)}$,
        $t'=\frac{a'_1}{r}+\sum_{i=2}^{\infty}\frac{a'_i}{(n'_i,m'_i)}$, $a_i$ and $a'_i$ as in  Proposition (3.4). Suppose that $a_i=a'_i, i=1,2,...,k-1$ and $a_k<a'_k$.\\
        If $|t'-t|<r^{-N} $ with $k<N$ then  $t=T_1+T_2+T_3,\ t'=T_1+T'_2+T'_3$ where\\
        $T_1=\displaystyle\frac{a_{1}}{r}+\sum_{{i=2} \atop
            {m_{i}+n_{i}=i}}^{k-1} \frac{a_{i}}{(n_{i},m_{i})},\
        T_3=\displaystyle\sum_{{i\geq N+1} \atop {m_{i}+n_{i}=i}}^{\infty}
        \frac{a_{i}}{(n_{i},m_{i})},\ T'_3=\displaystyle\sum_{{i\geq N+1}
            \atop {m_{i}+n_{i}=i}}^{\infty} \frac{a'_{i}}{(n'_{i},m'_{i})}$ and
        \\ \ \\
        $1)$ if $k$ is even and $a_{k}$ satisfies items $(3a)$ or $(3b)$ or
        $(3c)$ of
        Proposition \ref{l3} then\\ 
        $\begin{array}{l}
        T_2=\frac{a_{k}}{(n_{k},m_{k})}+\frac{r-1}{(n_{k}+1,m_{k})}+\frac{r-1}{(n_{k}+2,m_{k})}+\frac{r-3}{(n_{k}+2,m_{k}+1)}
        +\frac{r-1}{(n_{k}+3,m_{k}+1)}+\frac{r-3}{(n_{k}+3,m_{k}+2)}+...+\frac{a_{N}}{(n_{N},m_{N})},\\
        \ \\
        T'_2=\frac{a_{k}+1}{(n_{k},m_{k})}+\frac{0}{(n'_{k+1},m'_{k+1})}+\frac{0}{(n'_{k+2},m'_{k+2})}+\frac{0}{(n'_{k+3},m'_{k+3})}
        +\frac{0}{(n'_{k+4},m'_{k+4})}+\frac{0}{(n'_{k+5},m'_{k+5})}+...+\frac{0}{(n'_{N},m'_{N})},\end{array}$\\
        \ \\ Moreover we have $t=t'$ if and only if\\
        $t=T_1+\frac{a_{k}}{(n_{k},m_{k})}+\frac{r-1}{(n_{k}+1,m_{k})}+\sum_{i=0}^{\infty}\left(\frac{r-1}{(n_k+2+i,m_k+i)}+
        \frac{(r-3)}{(n_k+2+i,m_k+1+i)}\right)$ and \\ \ \\
        $t'=T_1+\frac{a_{k}+1}{(n_{k},m_{k})}+\sum_{i>k+1}^{\infty}\frac{0}{(n'_i,m'_i)}.$
        \ \\ 
        $2)$ if $k$ is odd and and $a_{k}$ satisfies item $(3c)$ of Proposition \ref{l3} then\\
        $\begin{array}{l}
        T_2=\frac{a_k}{(n_{k},m_{k})}+\frac{r-1}{(n_{k}+1,m_{k})}+\frac{r-3}{(n_{k}+1,m_{k}+1)}
        +\frac{r-1}{(n_{k}+2,m_{k}+1)}+\frac{r-3}{(n_{k}+2,m_{k}+2)}+...+\frac{a_{N}}{(n_{N},m_{N})},\\
        \ \\
        T'_2=\frac{a_{k}+1}{(n_{k},m_{k})}+\frac{0}{(n'_{k+1},m'_{k+1})}+\frac{0}{(n'_{k+2},m'_{k+2})}+\frac{0}{(n'_{k+3},m'_{k+3})}
        +\frac{0}{(n'_{k+4},m'_{k+4})}+...+\frac{0}{(n'_{N},m'_{N})},\end{array}$\\
        \ \\ Moreover we have $t=t'$ if and only if\\
        $t=T_1+\frac{a_{k}}{(n_{k},m_{k})}+\sum_{i=0}^{\infty}\left(\frac{r-1}{(n_k+1+i,m_k+i)}+
        \frac{(r-3)}{(n_k+1+i,m_k+1+i)}\right)$ and \\ \ \\
        $t'=T_1+\frac{a_{k}+1}{(n_{k},m_{k})}+\sum_{i>k+1}^{\infty}\frac{0}{(n'_i,m'_i)}.$
        \ \\  
        $3)$ if $k$ is odd and $a_{k}$ satisfies items $(4a)$ or $(4b)$ of
        Proposition \ref{l3} then\\
        $\begin{array}{l}
        T_2=\frac{a_k}{(n_{k},m_{k})}+\frac{r-3}{(n_{k},m_{k}+1)}+\frac{r-3}{(n_{k},m_{k}+2)}+\frac{r-1}{(n_{k}+1,m_{k}+2)}
        +\frac{r-3}{(n_{k}+1,m_{k}+3)}+\frac{r-1}{(n_{k}+2,m_{k}+3)}+...+\frac{a_{N}}{(n_{N},m_{N})},\\
        \ \\
        T'_2=\frac{a_{k}+1}{(n_{k},m_{k})}+\frac{0}{(n'_{k+1},m'_{k+1})}+\frac{0}{(n'_{k+2},m'_{k+2})}+\frac{0}{(n'_{k+3},m'_{k+3})}
        +\frac{0}{(n'_{k+4},m'_{k+4})}+\frac{0}{(n'_{k+5},m'_{k+5})}+...+\frac{0}{(n'_{N},m'_{N})},\end{array}$\\
        \ \\  Moreover we have $t=t'$ if and only if\\
        $t=T_1+\frac{a_{k}}{(n_{k},m_{k})}+\frac{r-3}{(n_{k},m_{k}+1)}+\sum_{i=0}^{\infty}\left(\frac{r-3}{(n_k+i,m_k+2+i)}+
        \frac{(r-1)}{(n_k+1+i,m_k+2+i)}\right)$ and \\ \ \\
        $t'=T_1+\frac{a_{k}+1}{(n_{k},m_{k})}+\sum_{i>k+1}^{\infty}\frac{0}{(n'_i,m'_i)}.$
        \ \\ 
        $4)$ if $k$ is even and $a_{k}$ satisfies item $(4c)$ of Proposition \ref{l3} then\\
        $\begin{array}{l}
        T_2=\frac{a_k}{(n_{k},m_{k})}+\frac{r-3}{(n_{k},m_{k}+1)}+\frac{r-1}{(n_{k}+1,m_{k}+1)}
        +\frac{r-3}{(n_{k}+1,m_{k}+2)}+\frac{r-1}{(n_{k}+2,m_{k}+2)}+...+\frac{a_{N}}{(n_{N},m_{N})},\\
        \ \\
        T'_2=\frac{a_{k}+1}{(n_{k},m_{k})}+\frac{0}{(n'_{k+1},m'_{k+1})}+\frac{0}{(n'_{k+2},m'_{k+2})}+\frac{0}{(n'_{k+3},m'_{k+3})}
        +\frac{0}{(n'_{k+4},m'_{k+4})}+...+\frac{0}{(n'_{N},m'_{N})},\end{array}$\\
        \ \\ Moreover we have $t=t'$ if and only if\\
        $t=T_1+\frac{a_{k}}{(n_{k},m_{k})}+\sum_{i=0}^{\infty}\left(\frac{r-3}{(n_k+i,m_k+1+i)}+
        \frac{(r-1)}{(n_k+1+i,m_k+1+i)}\right)$ and \\ \ \\
        $t'=T_1+\frac{a_{k}+1}{(n_{k},m_{k})}+\sum_{i>k+1}^{\infty}\frac{0}{(n'_i,m'_i)}.$
        \ \\ 
    \end{prop}

\begin{dem}
Take $t,t'\in[0,1]$ such that $|t'-t|<r^{-N}$,
$t=\displaystyle\frac{a_{1}}{r}+\sum_{{k=2} \atop
{m_{k}+n_{k}=k}}^{\infty} \displaystyle\frac{a_{k}}{(n_{k},m_{k})}$,
$t'=\displaystyle\frac{a'_{1}}{r}+\sum_{{k=2} \atop
{m'_{k}+n'_{k}=k}}^{\infty}
\displaystyle\frac{a'_{k}}{(n'_{k},m'_{k})}$, $a_i=a'_i,\ \forall\
i=1,2,...,k-1$ and $a_k<a'_k$.
Then \\
$$\displaystyle
t'-t=\frac{(a'_k-a_k)}{(n_{k},m_{k})}+\sum_{i>k}^{\infty}\left[\frac{a'_i}{(n'_i,m'_i)}-\frac{a_i}{(n_i,m_i)}\right]=\frac{(a'_k-a_k-1)}{(n_{k},m_{k})}+\frac{1}{(n_{k},m_{k})}+
\sum_{i>k}^{\infty}\left[\frac{a'_i}{(n'_i,m'_i)}-\frac{a_i}{(n_i,m_i)}\right],$$
and since $m_k+n_k=k$, $|t'-t|<r^{-N}$ then $a'_k-a_k-1=0$, that is,
$a'_k=a_k+1$.
\begin{enumerate}
\item Let $k$ be an even number and $a_{k}=0\ \mbox{or}\ 2n-1,\ n=1,...,a-1$.
Then $a_{k+1} \in \{0,1,...,r-1\}$ and we can write
$$\displaystyle\frac{1}{(n_k,m_k)}=\frac{r-1}{(n_k+1,m_k)}+
\displaystyle\sum_{i=0}^{\infty}\left(\frac{r-1}{(n_k+2+i,m_k+i)}+
\frac{(r-3)}{(n_k+2+i,m_k+1+i)}\right).$$ Therefore
$$\begin{array}{ll}
t'-t&\displaystyle=\frac{a'_{k+1}}{(n'_{k+1},m'_{k+1})}-\frac{a_{k+1}}{(n_k+1,m_k)}+
\frac{(r-1)}{(n_k+1,m_k)}+...=
\\ \ &\ \\ \ &
\displaystyle=\frac{a'_{k+1}}{(n'_{k+1},m'_{k+1})}+\frac{r-1-a_{k+1}}{(n_k+1,m_k)}+...,\end{array}$$
where $m'_{k+1}+n'_{k+1}=m_k+n_k+1=k+1$. As
$|t'-t|<r^{-N}<(r-2)^{-m}r^{-n},m+n=N\geq k+1$ then
$\frac{a'_{k+1}}{(n'_{k+1},m'_{k+1})}+\frac{r-1-a_{k+1}}{(n_k+1,m_k)}=0$
and it is possible only with $a'_{k+1}=0$ and $a_{k+1}=r-1$.\\
As $a_{k+1}=r-1$ and $k+1$ is an odd number, then $a_{k+2} \in
\{0,1,...,r-1\}$ and we have
$$\begin{array}{ll}
t'-t&\displaystyle=\frac{a'_{k+2}}{(n'_{k+2},m'_{k+2})}-\frac{a_{k+2}}{(n_k+2,m_k)}+
\frac{r-1}{(n_k+2,m_k)}+...=
\\ \ &\ \\ \ &\displaystyle=\frac{a'_{k+2}}{(n'_{k+2},m'_{k+2})}+\frac{r-1-a_{k+2}}{(n_k+2,m_k)}+....\end{array}$$
with $m'_{k+2}+n'_{k+2}=m_k+n_k+2=k+2$. Again we have $a'_{k+2}=0$ and $a_{k+2}=r-1$.\\
Now $a_{k+2}=r-1$ and $k+2$ is an even number. Then $a_{k+3} \in
\{0,1,...,(r-3)\}$ and
$$\begin{array}{ll}

t'-t&\displaystyle=\frac{a'_{k+3}}{(n'_{k+3},m'_{k+3})}-\frac{a_{k+3}}{(n_k+2,m_k+1)}+
\frac{(r-3)}{(n_k+2,m_k+1)}+...=\\
\ & \ \\
&\displaystyle=\frac{a'_{k+3}}{(n'_{k+3},m'_{k+3})}+\frac{(r-3)-a_{k+3}}{(n_k+2,m_k+1)}+....
\end{array}$$ with $m'_{k+3}+n'_{k+3}=m_k+n_k+3=k+3$. Therefore
$a'_{k+3}=0$ and $a_{k+3}=(r-3)$. Following this idea we have the
result.

\item To prove this part we use the same ideas of $(1)$ and the
equality
$$\displaystyle\frac{1}{(n_k,m_k)}=\displaystyle\sum_{i=0}^{\infty}\left(\frac{(r-1)}{(n_k+1+i,m_k+i)}+
\frac{(r-3)}{(n_k+1+i,m_k+1+i)}\right).$$

\item To prove this part we use the same ideas of $(1)$ and the
equality
$$\displaystyle\frac{1}{(n_k,m_k)}=\frac{(r-3)}{(n_k,m_k+1)}+\displaystyle\sum_{i=0}^{\infty}\left(\frac{(r-3)}{(n_k+i,m_k+2+i)}
+\frac{r-1}{(n_k+1+i,m_k+2+i)}\right).$$

\item To prove this part we use the same ideas of $(1)$ and the
equality
$$\displaystyle\frac{1}{(n_k,m_k)}=
\displaystyle\sum_{i=0}^{\infty}\left(\frac{(r-3)}{(n_k+i,m_k+1+i)}+\frac{r-1}{(n_k+1+i,m_k+1+i)}\right).$$
\hfill $\blacksquare$
\end{enumerate}
\end{dem}

Now we will give an explicit parametrization of
$\mathcal{B}_{\alpha-1}$ and hence for the boundary
$\partial\mathcal{R}_a$. Let $z$ be an element of
$\mathcal{B}_{\alpha-1}$. Using Proposition \ref{t2}, there exists a
sequence $(z_{n})_{n \geq 1}$ in $\mathcal{B}_{\alpha-1}$, such that
$$
z=g_{a_{1}} \circ g_{a_{2}} \circ \ldots \circ g_{a_{n}}
(z_{n}),\forall\  n\geq1.
$$
If $x$ is an element of $\mathcal{B'}_{\alpha-1}$, the sequence
$y_n=g_{a_{1}} \circ g_{a_{2}} \circ \ldots \circ g_{a_{n}}(x)$
converges to $z$ because the functions $g_{i}, i=0,1,...,r-1$ are
contractions.

Let $A=\{0,1,...,r-1\}$ be a subset of $\mathbb{N}$ and consider the
function
$$\begin{array}{ll}\psi:&A^{\mathbb{N}}\longrightarrow A^{\mathbb{N}}\\ \ &(a_i)\longmapsto \psi((a_i))=(b_i) \end{array}$$ given
by:\\
$\begin{array}{l} b_1=a_1;\\b_{2k}=r-1-a_{2k};\\b_{2k+1}=a_{2k+1}\
\mbox{if}\ a_{2k}\in\{0\}\cup \{2n-1:\ n=1...,a-1\};\\
b_{2k+1}=a_{2k+1}+2\ \mbox{if}\
a_{2k}\in\{2n:\ n=1,...,a-1\}.\end{array}$\\
Take $x_0\in \mathcal{B'}_{\alpha-1}$ and consider
$f:[0,1]\longrightarrow \mathcal{B}_{\alpha-1}$ defined as
follows:\\ if $t=\displaystyle\frac{a_{1}}{r}+\sum_{{k=2} \atop
    {m_{k}+n_{k}=k}}^{\infty}
\displaystyle\frac{a_{k}}{r^{n_{k}}(r-2)^{m_{k}}}, (a_i)\in
A^{\mathbb{N}}$, then $f(t)=\displaystyle\lim_{n\rightarrow
    \infty}g_{b_1}\circ g_{b_2}\circ...\circ
g_{b_n}(x_0)$ where $\psi((a_i))=(b_i)$.\\ \ \\
\begin{teo}\label{teofi} $f$ is a continuous, bijective function satisfying $f(0)=u=-1$
and $f(1)=v$.
\end{teo}
\dem\\
$(1)-$ $f$ is a well defined function.\\ We are going to use the
following notation:
$$g_{b_1}\circ g_{b_2}\circ...\circ g_{b_{k-1}}\circ g_{b_k}(x_0)=g_{b_1...b_k}(x_0).$$ According Lemma
\ref{l1} we have
$$\begin{array}{l}u=-1=g_{0(r-1)0(r-1)...}(x_0)=g_{\overline{0(r-1)}}(x_0).\\v=-(a-1)\alpha-\alpha^{-1}=g_{(r-1)0(r-1)0...}(x_0)=g_{\overline{(r-1)0}}(x_0).
\\w=-1-\alpha^3=g_{2(r-1)0(r-1)0...}(x_0)=g_{2\overline{(r-1)0}}(x_0).\end{array}$$
Taking $t,t'\in [0,1]$ such that $t=t'$. We have to prove that
$f(t)=f(t')$ and for this we use Proposition \ref{p37} and the
definition of $\psi$. We have to consider some cases.\\
-  $k$ be an even number and $a_{k}=0$ or $2n-1,\ n=1,...,a-1.$\\
Then $a_{k}+1=1$ or $2n$ and by Proposition \ref{p37} we
have\\
$t=\frac{a_1}{r}+...+\frac{a_{k-1}}{(n_{k-1},m_{k-1})}+\frac{a_k}{(n_k,m_k)}+\frac{r-1}{(n_k+1,m_k)}+
\sum_{i=0}^{\infty}\left(\frac{r-1}{(n_k+2+i,m_k+i)}+\frac{r-3}{(n_k+2+i,m_k+1+i)}\right)$\\
and\\
$t'=\frac{a_1}{r}+...+\frac{a_{k-1}}{(n_{k-1},m_{k-1})}+\frac{a_k+1}{(n_k,m_k)}+\sum_{i=k+1}^{\infty}\frac{0}{(n'_i,m'_i)}.$
Using the definition of $\psi$ we have\\ \ \\
$\begin{array}{l} f(t)=g_{b_1...b_{k-1}(r-1)\overline{(r-1)0}}(x_0)=g_{b_1...b_{k-1}(r-1)}(v),\\
f(t')=g_{b_1...b_{k-1}(r-2)\overline{0(r-1)}}(x_0)=g_{b_1...b_{k-1}(r-2)}(u)\end{array}$,\\
if
$a_k=0$ or\\
$f(t)=g_{b_1...b_{k-1}(r-2n)\overline{(r-1)0}}(x_0)=g_{b_1...b_{k-1}(r-2n)}(v)$\\
$f(t')=g_{b_1...b_{k-1}(r-2n-1)2\overline{(r-1)0}}(x_0)=g_{b_1...b_{k-1}(r-2n)}(w)$
if $a_k=2n-1,\ n=1,...,a-1$. \\ By Lemma \ref{l1} $f(t)=f(t').$\\
Using the same ideas we can prove the following cases (see \cite{JT}).\\
-  $k$ be an odd number and $a_{k}=2n-1,\ n=1,...,a-1$,\\
-  $k$ be an odd number and $a_{k}=0$ or $2n,\ n=1,...,a-2$,\\
-  $k$ be an even number and $a_{k}=2n,\ n=1,...,a-2$.\\
$(2)-$ $f$ is injective.\\Suppose that $f(t)=f(t')$. According to
Lemma \ref{l1} we have two possibilities:
\begin{itemize}
\item[$-$]  $f(t)=g_{b_1...b_{k-1}b_k}(u), b_k\in\{1,3,5,...,r-2\}$
and $f(t')=g_{b_1...b_{k-1}(b_k+1)}(v).$\\
Using the above notations we have
\\$f(t)=g_{b_1...b_{k-1}b_k\overline{0(r-1)}}(x_0)\ $ and
$f(t')=g_{b_1...b_{k-1}(b_k+1)\overline{(r-1)0}}(x_0).$ We need to
consider the following cases:
\begin{enumerate}
\item[$\bullet$] $k$ is an even number, $b_k\neq r-2$. In this case $b_k=r-1-a_k$ and then $a_k=r-1-b_k$ is an odd number.
By the definition of $\psi$ we have:\\
- $a_{k+1}=0$ because $b_{k+1}=0,\ a_k$ odd number,\\
- $a_{k+2}=0$ because $b_{k+2}=r-1,\ k+2$ even number.\\
Following this idea is easy to see that $a_i=0,\ \forall\ i\geq k+1$.\\
Therefore
$t=\frac{a_1}{r}+...+\frac{a_{k-1}}{(n_{k-1},m_{k-1})}+\frac{r-1-b_k}{(n_k,m_k)}+\sum_{i=k+1}^{\infty}\frac{0}{(n_i,m_i)}.$\\
We also have
$b'_k=b_k+1=r-1-a'_k$ and then $a'_k=r-2-b_k\neq 0$ is an even number. By the definition of $\psi$ and Proposition \ref{l3} we have:\\
-$a'_{k+1}=r-3,\ n'_{k+1}=n'_k,\ m'_{k+1}=m'_k+1$ because  $b'_{k+1}=r-1,\ a'_k$ even,\\
-$a'_{k+2}=r-1$ because $b'_{k+2}=0,\ k+2$ even ,\\
-$a_{k+3}=r-3,\ n'_{k+3}=n'_{k+2},\ m'_{k+3}=m'_{k+2}+1$.\\
Following this idea we have\\
$t'=\frac{a_1}{r}+...+\frac{a_{k-1}}{(n_{k-1},m_{k-1})}+\frac{r-2-b_k}{(n_k,m_k)}+
\sum_{i=0}^{\infty}\left(\frac{r-3}{(n'_k+i,m'_k+1+i)}+\frac{r-1}{(n'_k+1+i,m'_k+1+i)}\right)$
and then $t=t'$.
\end{enumerate}
Using the same ideas we can prove the following cases (see
\cite{JT}).
\begin{enumerate}
\item[$\bullet$] $k$ is an even number, $b_k=r-2$,
\item[$\bullet$]  $k$ is an odd number.

\end{enumerate}

\item[$-$] $f(t)=g_{b_1...b_{k-1}b_k}(-1-\alpha^3), b_k\in\{0,2,...,r-3\}$
and $f(t')=g_{b_1..b_{k-1}(b_k+1)}(v).$
Using the above notations we have
$f(t)=g_{b_1...b_{k-1}b_k2\overline{(r-1)0}}(x_0)\ $ and
$f(t')=g_{b_1...b_{k-1}(b_k+1)\overline{(r-1)0}}(x_0).$ We have to
consider the following cases:
\begin{enumerate}
\item[$\bullet$] $k$ is an even number. In this case $b_k=r-1-a_k$ and
$a_k=r-1-b_k$ is an even number too and we can prove that
$t=\frac{a_1}{r}+...+\frac{a_{k-1}}{(n_{k-1},m_{k-1})}+\frac{r-1-b_k}{(n_k,m_k)}+\sum_{i=k+1}^{\infty}\frac{0}{(n_i,m_i)}$\\
and\\
$t'=\frac{a_1}{r}+...+\frac{a_{k-1}}{(n_{k-1},m_{k-1})}+\frac{r-2-b_k}{(n_k,m_k)}+\frac{r-1}{(n'_k+1,m'_k)}+
\sum_{i=0}^{\infty}\left(\frac{r-1}{(n'_k+2+i,m'_k+i)}+\frac{r-3}{(n'_k+2+i,m'_k+1+i)}\right).$
Then $t=t'$


\item[$\bullet$] $k$ is an odd number. In this case $b_k=a_k$ or $b_k=a_k+2$ and then
$a_k=b_k$ or $a_k=b_k-2$. We can prove that
$t=\frac{a_1}{r}+...+\frac{a_{k-1}}{(n_{k-1},m_{k-1})}+\frac{a_k}{(n_k,m_k)}+\frac{r-3}{(n_k,m_k+1)}+
\sum_{i=0}^{\infty}\left(\frac{r-3}{(n_k+i,m_k+2+i)}+\frac{r-1}{(n_k+1+i,m_k+2+i)}\right)$\\ and\\
$t'=\frac{a_1}{r}+...+\frac{a_{k-1}}{(n_{k-1},m_{k-1})}+\frac{a_k+1}{(n_k,m_k)}+\sum_{i=k+1}^{\infty}\frac{0}{(n'_i,m'_i)}.$
Then $t=t'$.

\end{enumerate}
\end{itemize}
$(3)-$  $f$ is a continuous function. \\ Let us consider $t,t'\in
[0,1]$ 
, $|t'-t|<r^{-N}$ as in Proposition \ref{p37}. We have to consider
the following cases:
\begin{itemize}
\item [$-$]$t$ and $t'$ satisfying $(1)$ of Proposition \ref{p37}.\\
Here we have:
\begin{enumerate}
\item $f(t)=g_{b_1...b_{k-1}(r-1)(r-1)0(r-1)0...b_{N+1}...}(x_0)$ and
$f(t')=g_{b_1...b_{k-1}(r-2)0(r-1)0(r-1)...b'_{N+1}...}(x_0)$ if
$a_k=0$. Then\\
$\begin{array}{ll}|f(t)-f(t')|&=|g_{b_1b_2...b_{k-1}r-1}(z_1)-g_{b_1b_2...b_{k-1}r-2}(z_2)|\leq|\alpha|^{2(k-1)}|g_{r-1}(z_1)-g_{r-2}(z_2)|.\end{array}$\\
As $g_{r-2}(u)=g_{r-1}(v)$ then\\
$\begin{array}{ll}|f(t)-f(t')|&\leq
|\alpha|^{2(k-1)}\left(|g_{r-1}(z_1)-g_{r-1}(v)|+|g_{r-2}(z_2)-g_{r-2}(u)|\right)\leq\\
\ & \leq
|\alpha|^{2(k-1)}(|\alpha|^2+|\alpha|^3)diam(\mathcal{B}_{\alpha-1})=|\alpha|^{2k}(1+|\alpha|)diam(\mathcal{B}_{\alpha-1}),\end{array}$\\
where $diam(\mathcal{B}_{\alpha-1})$ is the diameter of $\mathcal{B}_{\alpha-1}$.\\
\item $f(t)=g_{b_1...b_{k-1}(r-2n)(r-1)0(r-1)0...b_{N+1}...}(x_0)$ and
$f(t')=g_{b_1...b_{k-1}(r-2n-1)2(r-1)0(r-1)0...b'_{N+1}...}(x_0)$ if
$a_k=2n-1$. Then\\
$\begin{array}{ll}|f(t)-f(t')|&=|g_{b_1b_2...b_{k-1}(r-2n)}(z_1)-g_{b_1b_2...b_{k-1}(r-2n-1)}(z_2)|\leq\\
\ &
\leq|\alpha|^{2(k-1)}|g_{r-2n}(z_1)-g_{r-2n-1}(z_2)|.\end{array}$\\
As $g_{r-2n-1}(w)=g_{r-2n}(v)$ then\\
$\begin{array}{ll}|f(t)-f(t')|&\leq
|\alpha|^{2(k-1)}\left(|g_{r-2n}(z_1)-g_{r-2n}(v)|+|g_{r-2n-1}(z_2)-g_{r-2n-1}(w)|\right)\leq\\
\ & \leq
|\alpha|^{2(k-1)}(|\alpha|^2+|\alpha|^3)diam(\mathcal{B}_{\alpha-1})=|\alpha|^{2k}(1+|\alpha|)diam(\mathcal{B}_{\alpha-1}).\end{array}$
\end{enumerate}
\end{itemize}

Using the same ideas we can prove the next cases (see \cite{JT}).
\begin{itemize}
    \item[$-$] $t$ and $t'$ satisfying item $(2)$ of Proposition \ref{p37}.

    \item[$-$] $t$ and $t'$ satisfying item $(3)$ of Proposition \ref{p37}.

    \item[$-$] $t$ and $t'$ satisfying item $(4)$ of Proposition \ref{p37}.\\

\end{itemize}

In all that cases we conclude that $f$ is a continuous
function.\hfill $\blacksquare$

Now we can finally prove Theorem \ref{princtheo}.

\dem\\ Let $\mathcal{Q}\subseteq \mathbb{R}^2$ be the set
$$\mathcal{Q}=\{(0,y),0\leq y\leq 1\}\cup \{(x,1),0\leq x\leq 1\}\cup \{(1,y),0\leq y\leq 1\}\cup \{(x,0),0\leq y\leq 1\}.$$Using Proposition \ref{t1}, Theorem \ref{teofi} and
Figure \ref{F4} we can prove that
$F:\mathcal{Q} \longrightarrow \partial \mathcal{R}_{a}$ given by\\ \ \\
$F(x,y)=\left\{\begin{array}{l}f(y),\mbox{if } (x,y)=(0,y),0\leq
y\leq 1\\
(f_2\circ f)(x),\mbox{if } (x,y)=(x,1),0\leq x\leq 1\\
(f_3\circ f)(y),\mbox{if } (x,y)=(1,y),0\leq y\leq 1\\
(f_1\circ f)(x),\mbox{if } (x,y)=(x,0),0\leq x\leq
1\end{array}\right.$ \ \\ is an homeomorphism . \hfill
$\blacksquare$

\begin{figure}[h]
    \centering \epsfig{file=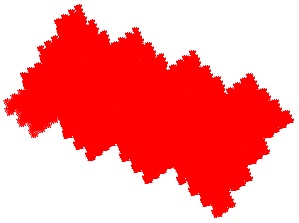,height=7cm,width=7.5cm} \caption{
        {\small $\mathcal{R}_4$ }}
\end{figure}

\begin{figure}[h]
    \centering \epsfig{file=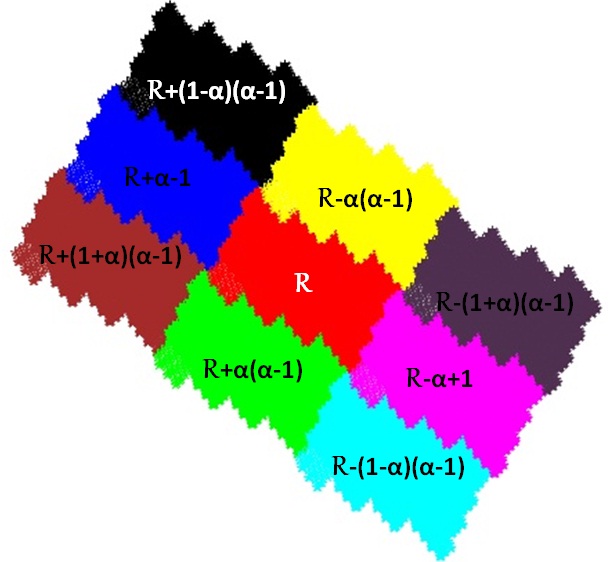,height=7cm,width=7.5cm}
    \caption{ {\small Tiling induced by $\mathcal{R}_4$}}
    \label{Figure1}
\end{figure}

\end{document}